\documentclass[12pt]{amsart}
\usepackage{amsfonts,amsmath,amssymb,amscd,amssymb}
\usepackage{enumerate}
\usepackage{enumitem}
\usepackage{comment}
\usepackage{hyperref}
\newtheorem{theorem}{Theorem}[section]
\newtheorem{lemma}[theorem]{Lemma}

\newtheorem{proposition}[theorem]{Proposition}

\newcommand{\Z}{\mbox{$\mathbb Z$}}
\newcommand{\Q}{\mbox{$\mathbb Q$}}
\newcommand{\N}{\mbox{$\mathbb N$}}     
\newcommand{\C}{\mbox{$\mathbb C$}}     

\begin{document}
	\title[]{Perfect powers with few digits and $S$-unit coefficients}
	
	\author[Darsana]{Darsana N}
	\address{Darsana N, Department of Mathematics, National Institute of Technology Calicut, 
		Kozhikode-673 601, India.}
	\email{darsana\_p230059ma@nitc.ac.in; darsanasrinarayan@gmail.com}
	
	\author[Nair]{Parvathi S Nair}%
	\address[Parvathi S Nair]{Department of Mathematics, National Institute of Technology Calicut, 
	Kozhikode-673 601, India.}
	\email[]{parvathisnair60@gmail.com, parvathi\_p220245ma@nitc.ac.in}
	
	\author[Rout]{S. S. Rout}
	\address{Sudhansu Sekhar Rout, Department of Mathematics, National Institute of Technology Calicut, 
		Kozhikode-673 601, India.}
	\email{sudhansu@nitc.ac.in; lbs.sudhansu@gmail.com}

	\dedicatory{}
\begin{abstract}
 In this paper, we study two questions concerning perfect powers in $S$-unit equations and sparse representations of perfect powers. Corvaja-Zannier \cite{corvaja2013finiteness} proved that there are only finitely many odd perfect powers in $\N$ having precisely four non-zero digits in their binary expansion. At first, we prove the finiteness of the set of solutions to the equation  
 \begin{equation*}
y^d=1+c_1g^{m_1}+c_2g^{m_2}+c_3g^{m_3}, \quad 0<m_1<m_2<m_3,
\end{equation*}
where $c_i$ are $S$-units for a suitable set $S$. The proof combines consequences of the Subspace Theorem concerning integer values of analytic series at $S$-unit points with Roth's theorem for number fields. As a second result, we extend the work of Bennett and Bugeaud \cite{bennettbugead2013perfect} by proving that, for every fixed integer base $g\geq 2$, there exists a constant $q_0(g)$ such that, for every prime $q>q_0(g)$, no perfect $q$-th power admits a representation with five non-zero digits in base $g$. The proof relies on explicit lower bounds for linear forms in both Archimedean and non-Archimedean logarithms.
\end{abstract}
\maketitle
\pagenumbering{arabic}
\pagestyle{headings}
\section{Introduction}
Diophantine equations involving perfect powers and sums of powers of a fixed integer constitute an important class of exponential Diophantine equations. It is known that precisely three numbers are currently known to have this property, namely
\begin{equation}\label{introeq1}
    3^4+3^3+3^2+3+1=11^2, 7^3+7^2+7+1=20^2, 18^2+18+1=7^3.
\end{equation}
It is widely conjectured that the corresponding Diophantine equation
\begin{equation}\label{introeq2}
    \frac{x^n-1}{x-1}=y^q
\end{equation}
in integers $x>1, y>1, n>2, q\geq2$ commonly known as the Nagell--Ljunggren equation, has no solutions other than \eqref{introeq1}. The case $q=2$ was completely settled by Ljunggren \cite{ljunggren1943noen}. Moreover, Ljunggren \cite{ljunggren1943noen} and Nagell \cite{nagell1920note} established that, if $3\mid n$ and $4\mid n$, then equation \eqref{introeq2} admits no solutions other than the three solutions given above. Shorey and Tijdeman \cite{ shorey1976new,exponential86} studied the finiteness of solutions under certain additional hypotheses. For further details on these finiteness results, as well as related developments, we refer the reader to their works and to the survey of Shorey \cite{ shorey2000exponential}. Equation \eqref{introeq2} has been completely resolved for several important classes of values of $x$. In particular, the equation has no solutions when $x$ is a square \cite{saradha1999equation, bugead1999}. It is also known that, when $x$ is a power of an integer belonging to $\{2,\ldots,10\}$, the only solutions are the two solutions listed above \cite{bugeaud1999integers}. Furthermore, under the hypothesis  stated in \cite{bugeaud2000diophantine,saradha1999equation}, the equation has been completely solved when $x$ is a power of a prime number. 

A related direction concerns perfect powers whose representations in a fixed base contain only a less number of non-zero digits. In the binary case, Szalay \cite{szalay2002equations} obtained a complete solution for odd squares represented by three binary digits, namely for equations of the form
\[
2^n\pm2^m\pm1=x^2.
\] This work has subsequently been extended and generalized in several directions. For instance, Luca \cite{luca2004diophantine} determined all positive integer solutions of
\[
x^2=p^a\pm p^b+1,
\]
where $p$ is an odd prime and $a>b>0$. The only solutions are
$(x,p,a,b)=(5,3,3,1),~(11,5,3,1)$. Corvaja and Zannier \cite{corvaja2013finiteness} proved that the Diophantine equation
\[
y^q=2^a+2^b+2^c+1
\]
admits only finitely many solutions in positive integers $(y,q,a,b,c)$. Their approach gives an effective upper bound for the exponent $q$, based on lower bounds for linear forms in logarithms. Moreover, for each fixed $q$, they established the finiteness of the corresponding solutions $(y,a,b,c)$. The proof combines Pad\'e approximations with an earlier result of Corvaja and Zannier \cite{corvaja2000diophantine} concerning $S$-unit points on analytic hypersurfaces, whose underlying arguments rely on Schmidt's Subspace Theorem. Bennett \cite{bennett2012perfect}, considered the analogous problem for ternary digits. He proved that if $y$ is coprime to $3$ and $y^2$ has at most three ternary digits, then $y\in\{1,5,8,13\}$ or $y=3^t+1$ for some nonnegative integer $t$. Bennett \cite{bennett2017exponential} also showed that, for any fixed prime $q$, the Diophantine equation
\[
6^a+2^b+1=y^q
\]
has only finitely many positive integer solutions $(y,a,b)$. Bennett, Bugeaud and Mignotte \cite{bennettbugead2012perfect, bennettbugead2013perfect}
further extended these results by completely determining the perfect powers whose base-$x$ representations, for $x=2$ or 3, are of the form $x^a+x^b+1$. Also, determined all perfect powers $y^q$, with $q\geq 4$, which can be represented in base $2$ as $y^q=2^a+2^b+2^c+1$. In their subsequent work \cite{bennettbugead2013perfect}, they generalized these investigations to perfect powers with few non-zero digits in an arbitrary fixed base . In particular, they proved that, for any fixed base $x\geq2$ and sufficiently large prime $q$, no perfect $q$-th power can be represented with three or four digits equal to $1$ in base $x$. More generally, for equations of the form
\[
x_1^a+x_2^b+1=y^q
\]
and
\[
x_1^a+x_2^b+x_3^c+1=y^q,
\]
where $x_1,x_2,x_3$ are fixed positive integers having a common prime divisor $p$ with $p\not\equiv1\bmod q$, they obtained an effective upper bound for $q$ in terms of $x_1,x_2,x_3$. Their results are obtained using refined lower bounds for linear forms in both Archimedean and non-Archimedean logarithms. This result has been sharpened in \cite{bennett2014perfect}. B{\'e}rczes et.al., \cite{berczes2025perfect} have been studied this finiteness results concerning perfect powers having two or three digits in their representation in a canonical number system of the equation of an algebraic number field.

In this direction, at first we extend a finiteness result of Corvaja and Zannier \cite{corvaja2005s} to powers of an arbitrary fixed integer with $S$-unit coefficients and further to the corresponding equations in $S$-units. Our main result is as follows.
\begin{theorem}\label{mainthm}
Suppose $d, g\geq 2$ are any given integers and $c_1, c_2, c_3$ are any fixed positive rational numbers which are $S$-units and $y\in \Z$. Assume that for any prime $p|g$, we have $(d,\varphi(p^n))=1$ for all $n\geq 1$, where $\varphi$ is the Euler totient function. Then 
\begin{equation}\label{eq4}
y^d=1+c_1g^{m_1}+c_2g^{m_2}+c_3g^{m_3}, \quad 0<m_1<m_2<m_3    
\end{equation}
has only finitely many solutions $(y, m_1, m_2, m_3)$ of \eqref{eq4} with $\gcd(g,c_i)=1$ for $i=1, 2, 3$ except in the following case:
\[
d=3,\qquad
c_1=3c,\quad c_2=3c^2,\quad c_3=c^3
\]
for some positive rational $S$-unit $c$, in which case the equation has the infinite family of solutions
\[
y=1+cg^m,\qquad
(m_1,m_2,m_3)=(m,2m,3m),
\]
whenever $1+cg^m\in\mathbb Z$.
\end{theorem}

Next, we consider a five-term extension of the above problem. We show that, under a suitable condition on the exponents, the prime exponent $q$ is effectively bounded. More precisely, our bound depends only on the parameter occurrs in this condition (see Remark 5.2 \cite{bennettbugead2013perfect}). 
\begin{theorem}\label{thm2}
Let $g \ge 2$ be a given integer, $m_1, m_2, m_3, m_4, y$ be positive integers, and $q$ be a prime number satisfying 
\[
y^q=1+g^{m_1}+g^{m_2}+g^{m_3}+g^{m_4}, \quad 0<m_1<m_2<m_3<m_4.
\]
 Assume that
 $q$ is coprime to $\varphi(g)$, and  there exists a real constant $\kappa>0$ such that $m_2 \ge \kappa m_3$.
Then the exponent $q$ is bounded by an effectively computable constant depending only on $\kappa$. In particular,
\[
q \le C(\kappa),
\]
for some constant $C(\kappa)>0$ depending only on $\kappa$.
\end{theorem}

\section{Preliminaries}
Let $K$ be a number field and let $S$ be a finite set of places of $K$ containing all the archimedean places. We denote $M_K$ for the set of places of $K$. For an $x\in K^{\times}$, the absolute logarithmic Weil height is defined as 
\[h(x)=\sum_{\nu \in M_k}\max\{0, \log|x|_{\nu}\},\]
where $|\cdot|_{\nu}$ denotes the normalized absolute value with respect to $K$ for a place $\nu$ of $K$. 
In particular, if $x$ is a nonzero rational with $x = m/n$ for $m$ and $n$ coprime integers, we define the logarithmic height of $x$ as $h(x) = \max\{\log|m|,\log|n|, 1\}$.

For $x\in K^{\times}$, the $S$-height is defined as 
\[h_{S}(x)=\sum_{\nu \notin S}\max\{0, \log|x|_{\nu}\}.\]

For a vector ${\bf x}=(x_1,\ldots, x_n)\in K^{n}\backslash \{0\}$, the projective height is defined as
\[h({\bf x})=\sum_{\nu \in M_K}\log \max\{|x_1|_{\nu},\ldots, |x_n|_{\nu}\}.\]
Also,
\[\hat{h}({\bf x})=\sum_{1\leq i\leq n}h(x_i).\]
In the same sense, we can define $\hat{h}_S({\bf x})$ as sum over all $h_S(x_i)$.
For an absolute value $\nu \in M_K$, we denote by $\C_{\nu}$ a completion of an algebraic closure of $K_{\nu}$.

We recall the structure of an algebraic subgroup of $G_m^{n}$ in the sense of \cite{corvaja2005s}. If $H$ is an algebraic subgroup of $G_m^{n}$, then $H$ is defined by equations of the form $X_1^{a_1}\cdots X_n^{a_n}=1$, where $(a_1,\ldots, a_n)$ varies over a lattice $\Lambda \subset \Z^n$. Also, $\dim H+rank\Lambda=n$. $H$ is connected if $\Z^n/\Lambda$ is torsion-free. When $H$ is connected, we call $H$ a torus. By a coset, we mean a translate of a torus by a point.

When $H$ is a torus, we have an algebraic group isomorphism between $H$ and $G_m^{\dim H}$ given by the parametrisation
\begin{equation}\label{eqtorus}
X_i={\bf t}^{{\bf u}_i}, ~~ i=1,\ldots, n, ~~t=(t_1,\ldots, t_s),    
\end{equation}
where $s=\dim H$ and ${\bf u}_i\in \Z^s$. We say an algebraic subgroup $H \subset G_m^n$ passes through $0$ if the Zariski closure of $H$ in $\mathbb{A}^n$ contains $(0,\ldots, 0)$. Next proposition gives equivalent conditions for a torus to pass through $0$.

\begin{proposition}[\cite{corvaja2005s}, Proposition 1]\label{proptorus}
For a torus $H$, the following are equivalent. 
\begin{enumerate}[label=\textnormal{(\roman*)}]
\item $H$ passes through $0$.
\item The lattice $\Lambda$ does not contain any non-zero vector with all its coordinates $\geq 0$.
\item There exists a parametrisation as in \eqref{eqtorus} with the vectors ${\bf u}_i \in \Z^s$ have all their coordinates $>0$.
\item There exists a point $(x_1,\ldots, x_n)\in H\cap \C_{\nu}^n$ with $|x_i|_{\nu}<1$ for $i=1,\ldots, n$.
\item There exists a sequence in $H\cap \C_{\nu}^n$ converging to $(0,\ldots, 0)$ in the $\nu$-adic topology.
\end{enumerate}
\end{proposition}

We adopt the following theorem from \cite{corvaja2005s}.
\begin{theorem}[\cite{corvaja2005s}, Theorem 1]\label{thmseries}
Consider the power series $g({\bf X})=\sum_{\bf i}a_{\bf i}{\bf X}^{\bf i}$ with algebraic coefficients in $\C_{\nu}$ converging in an neighborhood of the origin in $\C_{\nu}^n$. Let $S$ be a finite set of places of $K$ containing all the archimedean places. Let ${\bf x}_h=(x_{h1},\ldots, x_{hn})$ for $h=1, 2, \ldots$, be a sequence in $({K^{\times}})^n$, tending to zero in $K_{\nu}^{n}$ and such that $g({\bf x}_h)$ is defined and belongs to $K$. Assume that 
\begin{enumerate}[label=\textnormal{(\roman*)}]
\item for $i=1,\ldots, n$, we have $h_S(x_{hi})+h_S(x_{hi}^{-1})=o(h(x_{hi}))$ as $h\rightarrow \infty$.
\item $\hat{h}({\bf x}_h) =O(-\log \max_i|x_{hi}|_{\nu})$.
\item $h_S(g({\bf x}_h))=o(h({\bf x}_h))$.
\item $h(g({\bf x}_h))=O(h({\bf x}_h))$.
\end{enumerate}
Then there exists a finite number of cosets ${\bf u}_1H_1,\ldots, {\bf u}_rH_r\subset G_m^{n}$ such that $\{{\bf x}_h\}_{h\in N}\subset \cup_{i=1}^{r}{\bf u}_i H_i$ and such that, for $i =1,\ldots, r$, the restriction of $g({\bf X})$ to ${\bf u}_iH_i$ coincides with a
polynomial in $K[{\bf X}]$.
\end{theorem}
The following theorem was conjectured by Lang and proved by Laurent.
\begin{theorem}[\cite{bombieri2006heights}, Theorem 7.4.7]\label{thmlang}
Let $\Gamma$ be a finitely generated subgroup of $(\overline{\Q}^{\times})^n$ and let $\Sigma$ be any subset of $\Gamma$. Then the Zariski closure of $\Sigma$ in $G_m^{n}$ is a finite union of translates of algebraic subgroups of $G_m^{n}$.
\end{theorem}
The lemma below is adopted from \cite{corvaja2013finiteness}.

\begin{lemma}[\cite{corvaja2013finiteness}, Lemma 2.1]\label{lemmaseries}
For $d\geq 2$, $\zeta_1, \zeta_2, \zeta_3 \in \C^{\times}$ and for positive integers $l_1\leq l_2\leq l_3$, the polynomial $1+\zeta_1T^{l_1}+\zeta_2T^{l_2}+\zeta_3T^{l_3}$ is a perfect $d$-th power of a non-constant complex polynomial precisely in the following cases.
\begin{enumerate}
\item $d=2, l_1=l_2<l_3$ and $4\zeta_3=(\zeta_1+\zeta_2)^2$.
\item $d=2, l_1<l_2=l_3$ and $4(\zeta_2+\zeta_3)=\zeta_1^2$.
\item $d=2, l_2=3l_1, l_3=4l_1$ and $8\zeta_2=-\zeta_1^3, 64\zeta_3=\zeta_1^4$.
\item $d=3, l_2=2l_1, l_3=3l_1$ and $3\zeta_2=\zeta_1^2, 27\zeta_3=\zeta_1^3$.
\end{enumerate}
\end{lemma}

As usual, $\Gamma(z)$ shall denote Euler's Gamma-function, while
$(1+z)^{\frac{1}{d}}$ shall mean the binomial series
\[
\sum_{n=0}^{\infty} \binom{\frac{1}{d}}{n} z^n \in \mathbb{Q}[[z]],
\]
interpreted for the moment in the formal sense; also, we shall adopt the
following standard notation for a hypergeometric function:
\begin{equation}\label{prelem2eq1}
F(a,b,c,z)=1+\frac{a\cdot b}{1\cdot c}z
+\frac{a(a+1)\cdot b(b+1)}{1\cdot 2\cdot c(c+1)}z^2+\cdots.
\end{equation}

\begin{lemma}[\cite{corvaja2013finiteness}, Lemma 2.4]\label{prelem2}
Let $r,s$ be positive integers and set $t:=r+s$. Also, put
\[
G(z):=
\frac{\Gamma\left(\frac{1}{d}+1\right)\Gamma(1+t)}
{\Gamma\left(\frac{1}{d}+r+1\right)\Gamma(s+1)}
F\left(\frac{1}{d}-r,-s,-t,z\right),
\]
\[
H(z):=
F\left(\frac{1}{d}-s,-r,\frac{1}{d}+1-z\right),
\]
and
\[
E(z):=
(-1)^{s+1}
\frac{\Gamma\left(\frac{1}{d}+1\right)r!}
{\Gamma\left(\frac{1}{d}-s\right)t!}
F\left(-\frac{1}{d}+s+1,s+1,t+2,z\right).
\]
Then $G(z)$ and $H(z)$ are polynomials with rational coefficients, of
degree $s$, resp.\ $r$ and
\begin{equation}\label{prelem2eq2}
G(z)-(1-z)^{\frac{1}{d}}H(z)=z^{t+1}E(z).
\end{equation}
Finally, there is a positive integer $A$ such that $AG(z)$ and $AH(z)$
have integer coefficients bounded in absolute value by
$B:=(2d)^{4t}$. More precisely, we may take
\[
A=d^{4t}\binom{\frac{1}{d}+r}{r}.
\]
\end{lemma}

Theorem \ref{Roth} is the general form of Roth's theorem over number fields.
\begin{theorem}[\cite{bombieri2006heights}, Theorem 6.2.3]\label{Roth}
Let $K$ be a number field with a finite set $S$ of places of $K$. For $\nu \in S$, let $\alpha_{\nu}\in K_{\nu}$ be $K$-algebraic. Let $\kappa\geq 2$. Then there are only finitely many $\beta \in K$ such that 
\[\prod_{\nu \in S}\min\{1, |\beta-\alpha_{\nu}|_{\nu}\}\leq H(\beta)^{-\kappa}.\]
\end{theorem}

We need the following special version of a corollary obtained in \cite{laurent2008linear} to prove Theorem \ref{thm2}.
\begin{proposition}\label{preprop1}
    Let $\alpha_1$ and $\alpha_2$ be multiplicatively independent positive rational numbers, and $b_1$ and $b_2$ be positive integers. Define
    \[
    \Lambda=|b_2\log \alpha_2-b_1\log \alpha_1|.
    \] Then
    \[
    \log \Lambda\geq -25.2\left( \max\{\log b^{'}+0.38,10\}\right)^2h(\alpha_1)h(\alpha_2),
    \] where
    \[
    b^{'}=\frac{b_1}{h(\alpha_2)}+\frac{b_2}{h(\alpha_1)}.
    \]
\end{proposition}
We now turn our attention to non-Archimedean result. Suppose we are given
two multiplicatively independent positive rational numbers $\alpha_1$ and $\alpha_2$. We define $g$ to be the smallest positive integer such that both $\nu_p(\alpha_1^g-1)$ and $\nu_p(\alpha_2^g-1)$ are positive.
Note that if $m=p_1^{j_1}\cdots p_k^{j_k}$ where $p_i'$s are distinct primes and $j_i\in \N,$ we define, for a nonzero integer $x$,
\[
\nu_m(x)=\min_{1\leq i\leq k}\left[\frac{\nu_{p_i}(x)}{j_i}\right].
\]
\begin{proposition}[\cite{bugeaud2002linear}, Theorem 3]\label{preprop2}
    Let $\alpha_1$ and $\alpha_2$ be positive rational numbers with $\alpha_1\neq 1$, $b_1$ and $b_2$ be positive integers and set
    \[
    \Lambda=\alpha_2^{b_2}-\alpha_1^{b_1}.
    \] For any set of distinct primes $p_1,\ldots,p_k$ and positive integers $j_1,\ldots,j_k$, we set $m=p_1^{j_1}\cdots p_k^{j_k}$ and suppose that there exists a positive integer $g$ such that for each $i$, we have 
    \[
    \nu_{p_i}(\alpha_1^g-1)\geq j_i
 \text{ and  }   \nu_{p_i}(\alpha_2^g-1)\geq 1, \text{ if } p_i\geq 2, \] and also
 \[
 \nu_{p_i}(\alpha_1^g-1)\geq 2 \text{ and } \nu_{p_i}(\alpha_2^g-1)\geq 2, \text{ if } p_i=2.
 \] Then, if $m,b_1$ and $b_2$ are relatively prime, we have
 \[
 \nu_m(\Lambda)\leq \frac{66.8g}{(\log m)^4}\left( \max\{\log b^{'}+\log (\log m)+0.64, 4 \log m\}\right)^2(\log A_1)(\log A_2),
 \] where 
 \[
 b^{'}=\frac{b_1}{\log A_2}+\frac{b_2}{\log A_1} \text{ and } \log A_i\geq \max\{h(\alpha_i), \log m\}.
 \]
\end{proposition}

\begin{lemma}[\cite{bennettbugead2013perfect}, Lemma 4.1]\label{prelem1}
    Let $y\geq 2$ be an integer and $q$ be a prime. Assume that $y^q\equiv 1\pmod q$. Then $y\equiv 1\pmod q$ and $y^q\equiv 1\pmod {q^2}$. Moreover, if in addition $q\geq 3$, then $\nu_q(y-1)=\nu_q(y^q-1)-1$.
\end{lemma}

\section{Proof of Theorem
\ref{mainthm}.}
We fix the set $S$ to be a finite set containing all the archimedean places and places above the prime factors of the numbers $c_1, c_2, c_3$ and $g$. 
Without loss of generality, we can assume that $d$ is a prime number. We split the proof into four cases as follows.
\subsection{Case-1}
There exists an infinite subsequence of solutions such that $m_2\leq\kappa_1 m_3$, where $9/10<\kappa_1<1$.

Dividing \eqref{eq4} by $c_3g^{m_3}$ to obtain 
\[\frac{1}{c_3}(g^{m_3})^{-1}y^d=1+\frac{c_2}{c_3}g^{m_2-m_3}+\frac{c_1}{c_3}g^{m_1-m_3}+\frac{1}{c_3}g^{-m_3}.\]
Put $n_1=m_3-m_2, n_2=m_3-m_1, n_3=m_3$, $x_1=c_2/c_3 g^{-n_1}, x_2=c_1/c_3 g^{-n_2}$, and $x_3=1/c_3 g^{-n_3}$, then
\begin{equation}
\frac{1}{c_3}(g^{m_3})^{-1}y^d=1+x_1+x_2+x_3.    
\end{equation}
The series defined by
\begin{align}\label{eqseries}
\begin{split}
F(X_1,X_2,X_3)&= \\&=1+\frac{1}{d}(X_1+X_2+X_3)+\frac{(1-d)}{2d^2}(X_1+X_2+X_3)^2+\cdots.    
\end{split}
\end{align}
is certainly converging absolutely to a continuous function for $X_i \in \C$ with $|X_i|\leq 1/4$ $(i=1,2,3)$. Choose $\kappa_2$ such that $0<\kappa_2<\kappa_1$ and $\kappa_1+\kappa_2=1$. We observe that $n_3>n_2>n_1\geq \kappa_2 m_3$ and, by our assumption, $m_2\leq \kappa_1m_3$. Hence, for large values of $m_3$, we have $|x_i|\leq 1/4$ for $i=1,2,3$. Thus, the series $F$ converges at $(x_1, x_2, x_3)$ and we may write $z=F(x_1, x_2, x_3)$ and 
\begin{equation}\label{eqzd}
z^d=1+x_1+x_2+x_3=1+\frac{c_2}{c_3}g^{m_2-m_3}+\frac{c_1}{c_3}g^{m_1-m_3}+\frac{1}{c_3}g^{-m_3}.    
\end{equation}
We choose $K$ as the splitting field of $X^d-gc_3$ over $\Q$, so that the sequence of solutions to \eqref{eqzd} is defined over $K$. Note that $z=({c_3g^{m_3}})^{\frac{-1}{d}}y$. We let $S$ be a finite set of places of $K$ that includes all the archimedean places, all places above the prime divisors of $g$, and all places above prime divisors of $c_1, c_2$ and $c_3$. Hence, we see that $x_1, x_2, x_3$ are $S$-units and $z$ is an $S$-integer. The complex absolute value induces an absolute value on $\Q(z)$, and we extend this to an infinite place of $K$. So we may embed $K$ in $\C$ by means of $\nu$, and still, we can write $z=F(x_1, x_2, x_3)$ with respect to the $\nu$-adic convergence.

We have that
$h_S(x_i)=h_S(x_i^{-1})=0$ as $x_i$ are $S$-units for $i=1,2,3$. Hence, condition $(i)$ of Theorem \ref{thmseries} holds. Now to verify
$\hat{h}(x_1,x_2,x_3)=O\left( -\log \max_{i}|x_i|_{\nu}\right)$, where $\hat{h}(x_1,x_2,x_3)=\sum_{i=1}^3h(x_i)$, we first find $h(x_1)$.
\begin{align*}
    h(x_1)=h\left(\frac{c_2}{c_3}g^{-n_1}\right)\leq h\left(\frac{c_2}{c_3}\right)+n_1h\left(g\right)\ll m_3.
\end{align*} Similarly, we get $h(x_i)\ll m_3$ for $i=2,3$. Therefore, 
\begin{equation}\label{propeq1}
    \sum_{i=1}^3h(x_i)\ll 3m_3.
\end{equation} 
Next, we will show that $m_3\ll  -\log \max_{i}|x_i|_{\nu} $. We have
\begin{align*}
    |x_1|_{\nu}=\left|\frac{c_2}{c_3}g^{-n_1}\right|_{\nu}=\left|\frac{c_2}{c_3}\right|_{\nu}|g^{m_2-m_3}|_{\nu}\leq \left|\frac{c_2}{c_3}\right|_{\nu}|g^{-\kappa_2m_3}|_{\nu}.
\end{align*}
In a similar way, we obtain a bound for $|x_2|_{\nu}$ and $|x_3|_{\nu}$ in terms of the power of $-m_3$. Therefore,
\[
\max_{i}|x_i|_{\nu}\leq \max\left\{\frac{c_2}{c_3},\frac{c_1}{c_3},\frac{1}{c_3}\right\}g^{-\kappa_2m_3},
\]
which implies 
\[
-\log \left(\max_{i}|x_i|_{\nu}\right)\geq -\log A+\kappa_2m_3\log g,
\] where $A:=\max\left\{\frac{c_2}{c_3},\frac{c_1}{c_3},\frac{1}{c_3}\right\}$. Thus, we obtain
\begin{equation}\label{propeq2}
    m_3\ll-\log \left(\max_{i}|x_i|_{\nu}\right),
\end{equation}
 since $\log A\leq -\log \left(\max_{i}|x_i|_{\nu}\right)$. From \eqref{propeq1} and \eqref{propeq2}, we conclude that $\hat{h}(x_1,x_2,x_3)=O\left( -\log \max_{i}|x_i|_{\nu}\right)$.
Since $z$ is an $S$-integer, $h_S(F(x_1, x_2, x_3))=h_S(z)=0$. This gives that $h_S(F(x_1,x_2,x_3))=o\left(h(x_1,x_2,x_3)\right)$. Further,
$h(F(x_1,x_2,x_3))=h(z)\ll m_3$,

Now we apply Theorem \ref{thmseries} to the series $F$ and to the sequence of $S$-unit points $(x_1, x_2, x_3)$. Then we obtain a finite number of cosets $u_1H_1,\ldots, u_rH_r$ of $G_m^{3}$, with $u_i\in G_m^{3}(K)$ with $H_i$ connected algebraic subgroup of $G_m^{3}$ such that $(x_1, x_2, x_3)\in \cup_{i=1}^{r}{u}_i H_i$ and the restriction of $F(X_1, X_2, X_3)$ to $u_iH_i$ coincides with a polynomial.

When considering an infinite sequence of solutions, we may assume that there exists a coset that contains all our points. Also, our sequence consists of $S$-units, so we can apply Theorem \ref{thmlang} to our sequence. So, we may suppose that the sequence is Zariski-dense in our coset. Denote such coset by $uH$ with $u=(\zeta_1, \zeta_2, \zeta_3)\in G_m^{3}$ being any of our points. This coset is not a single point, so $\dim H=:s>0$. We have that our points converge $\nu$-adically to the origin. Hence, by the Proposition \ref{proptorus}, we get a parametrization of $uH$ given by 
\begin{equation}\label{eqX}
X_i={\zeta}_i T_1^{a_{i1}}\cdots T_s^{a_{is}}    
\end{equation} 
for $i=1, 2, 3$ and $a_{ij}\geq 0$ for all $i, j$. Also, we have that the restriction of $F(X_1, X_2, X_3)$ to $u_iH_i$ coincides with a polynomial. Hence, $F(X_1, X_2, X_3)$ is a polynomial in $T_1,\ldots, T_s$. Since our sequence is Zariski dense in the coset, for $X_i$ in \eqref{eqX} we have that \[F^d(X_1, X_2, X_3)=1+X_1+X_2+X_3.\]

It is clear that $X_i$ cannot be identically $1$ on $H$. Now for the non-zero vectors $(a_{i1},\ldots, a_{is})\in \N^s$ for $i=1, 2, 3$, we can find positive integers $b_1,\ldots, b_s$ such that the quantity $l_i=\sum_{k=1}^sb_ka_{ik}$ are all positive and such that for $i, j \in \{1, 2, 3\}$, $l_i=l_j$ if and only if $(a_{i1},\ldots, a_{is})=(a_{j1},\ldots, a_{js})$. By substituting $T^{b_j}$ for $T_j$, we obtain a non-constant polynomial $P(T)$ such that 
\begin{equation}\label{eqpoly}
P^d(T)=1+\zeta_1T^{l_1}+\zeta_2T^{l_2}+\zeta_3T^{l_3}.    
\end{equation} By renumbering we may assume that $0<l_1\leq l_2\leq l_3$. Hence, we can use Lemma \ref{lemmaseries} on the polynomial \eqref{eqpoly}. Let us write the elements of $S$ explicitly as $\{p_1,p_2,\ldots,p_r\}$. Then we may write
\[
\zeta_1=p_1^{a_1}\cdots p_r^{a_r}, ~ \zeta_2=p_1^{b_1}\cdots p_r^{b_r},~ \zeta_3=p_1^{c_1}\cdots p_r^{c_r},
\] where $a_i,b_i,c_i\in\mathbb{Z}$.
Suppose we fall in case(A) of Lemma \ref{lemmaseries}. Then
\[
4\zeta_3=4p_1^{c_1}\cdots p_r^{c_r}=(p_1^{a_1}\cdots p_r^{a_r}+p_1^{b_1}\cdots p_r^{b_r})^2.
\]
Let $d_i=\min\{a_i,b_i\}$. Then
\begin{equation}\label{eqzeta}
4\zeta_3=4p_1^{c_1}\cdots p_r^{c_r}=p_1^{2d_1}\cdots p_r^{2d_r}(p_1^{a_1-d_1}\cdots p_r^{a_r-d_r}+p_1^{b_1-a_1}\cdots p_r^{b_r-a_r})^2. 
\end{equation}
Now for any $p_k\in S$, $p_k\nmid (p_1^{a_1-d_1}\cdots p_r^{a_r-d_r}+p_1^{b_1-a_1}\cdots p_r^{b_r-a_r})$. Thus, we can write
\[
p_1^{a_1-d_1}\cdots p_r^{a_r-d_r}+p_1^{b_1-a_1}\cdots p_r^{b_r-a_r}=q_1^{e_1}\cdots q_s^{e_s}\geq 4, 
\]  where $q_j\notin S$ for every $1\leq j\leq s$.
Now we consider the case when $S$ contains $2$. Then from \eqref{eqzeta}, $q_j$ must be equal to $p_i$ for some $i$ and $j$, which is not possible. On the other hand, if $S$ does not contain $2$, then one of the $q_j$ must be equal to 2. Without loss of generality, assume that $q_1=2$. Then \eqref{eqzeta} becomes
\[
\zeta_3=p_1^{2d_1}\cdots p_r^{2d_r}(q_1^{e_1-1}\cdots q_s^{e_s})^2.
\]
But this will contradict the fact that $\zeta_3$ is an $S$-unit. Similarly, for the case(B). Further, we cannot fall in case(C) since $\zeta_1,\zeta_2$ are positive.
Finally, when we are in case(D), then \eqref{eqpoly} becomes
\[
P^3(T)=1+\zeta_1T^{\ell_1}+\frac{\zeta_1^2}{3}T^{2\ell_1}+\frac{\zeta_1^3}{27}T^{3\ell_1}=\left(1+\frac{\zeta_1}{3}T^{\ell_1}\right)^3.
\]
Thus, for $d=3$, if we choose $c_1=3c, c_2=3c^2$ and $c_3=c^3$, for some appropriate positive rational $S$-unit $c$, we get an infinite set of integer solutions 
\[
y=1+cg^m,\qquad
(m_1,m_2,m_3)=(m,2m,3m),\quad m\in \N
\]
to the equation \eqref{eq4}.

\subsection{Case-2} There exists an infinite subsequence of solutions such that $m_1\geq \kappa_2m_3$. 
We have 
\[y^d=1+x_1+x_2+x_3,\]
with $x_i=c_ig^{m_i}$ for $i=1, 2, 3$.
Let $g=p_1^{k_1}\cdots p_r^{k_r}$ be the prime factorisation of $g$. Hence, $x_i=c_i(p_1^{k_1}\cdots p_r^{k_r})^{m_i}$. By assumption, $(c_i, g)=1$ and $m_1\geq \kappa_2m_3$. Hence, for any $p_j, j=1,\ldots, r$,
\[|x_i|_{p_j}=|c_i(p_1^{k_1}\cdots p_r^{k_r})^{m_i}|_{p_j}=p_j^{-k_jm_i}\leq p_j^{-\kappa_2k_j m_3}.\]
For large values of $m_3$, we have $|x_i|_{p_j}\rightarrow 0$. Hence, the series in \eqref{eqseries} converges in $\Q_{p_j}$ at $(x_1, x_2, x_3)$. We let $z=F(x_1, x_2, x_3)$. Thus
\[z^d=1+x_1+x_2+x_3\]
and $z=y\theta$ for some $d$-th root of unity $\theta$.
Corresponding to each $p_j$, there exists a place in $\Q(z)$. We extend this place to a place $\nu$ of $K$ lying above $p_j$. We apply Theorem \ref{thmseries} with this $\nu$. Since $x_i$ and $x_i^{-1} $ are $S$-units, we have $h_S(x_i)+h_S(x_i^{-1})=0$. Thus we get $h_S(x_i)+h_S(x_i^{-1})=o(h(x_i))$. We have $h(x_i)\leq h(c_i)+h(g^{m_i})$. Thus
\[
\hat{h}((x_1,x_2,x_3))=\sum_{i=1}^3h(x_i)\ll 3m_3.
\] Also, $\max\{|x_i|_\nu\}\leq p_j^{-\kappa_2k_jm_3}$ implies that $m_3\ll -\log (\max\{|x_i|_\nu\})$ which gives the condition (ii) of Theorem \ref{thmseries}. For any $\mu\notin S $, $|z^d|_\mu\le \max\{|x_1|_{\mu}, |x_2|_{\mu}, |x_3|_{\mu}\}\leq 1$. So $z$ is an $S$-integer which shows that $h_S(z)=0$ and hence condition (iii) of Theorem \ref{thmseries} holds. Similarly we get $h(F(x_1,x_2,x_3))=O(h(x_1,x_2,x_3))$. We follow the same arguments as in the first case and apply Lemma \ref{lemmaseries} and obtain contradiction in each case.

Now onward, we will restrict our search for solutions satisfying 
\begin{equation}\label{eqrestr}
0<m_1<\kappa_2m_3 \text{ and } \kappa_1m_3<m_2<m_3.
\end{equation}
\subsection{Case-3} There exists an infinite subsequence of solutions such that $m_1$ is bounded. 

Suppose that $m_1=b$ on an infinite subsequence, where $b$ is a constant. From \eqref{eq4}, we have 
\begin{equation}\label{eqdelta}
y^d-1-c_1g^{m_1}=c_2g^{m_2}+c_3g^{m_3}.    
\end{equation} 
Substituting $\Delta=1+c_1g^{m_1}$ in \eqref{eqdelta}, we get
\begin{equation}\label{eq10}
y^d-\Delta=g^{m_2}(c_2+c_3g^{m_3-m_2}).    
\end{equation}
We get that there exists a fixed $d$-th root of unity $\theta$ of $\Delta$ such that $\Delta=\theta^d$. 
Let the prime factorisation of $g$ be 
\[g=p_1^{a_1}\cdots p_k^{a_k}, a_i\in \Z \text{ for } i=1,\ldots, k.\]
Since $c_i$ and $g$ are coprime for $i=1, 2, 3$ then for any prime $p_j$ for $j=1,\ldots, k$ equation \eqref{eq10} implies that 
\[|y^d-\Delta|_{p_j}=p_j^{-a_jm_2}.\]
Using 
\[|y^d-\theta^d|_{p_j}=|y-\theta|_{p_j}|y^{d-1}+y^{d-2}\theta+\cdots+\theta^{d-1}|_{p_j}\] and the ultra-metric triangle inequality, we get that 
$|y-\theta|_p\leq p^{-am_2}$ for some constant $a\in \Z$. Now, for any prime $q$ dividing any of the $c_i$ for $i=1,2,3$, we get $|y-\theta|_q\leq q^{-b}\ll 1$ for some constant $b\in \Z$. We use Theorem \ref{Roth} for the finite set $S$ that contains all the archimedean places and places above the prime factors of the numbers $c_1, c_2, c_3$, and $g$, to derive a contradiction.              

For the place above the archimedean place, we fix $\alpha_{\infty}=\infty$. As the distance from $y$ to $\infty$ is meaningless, we set $1/y=y-\infty$ and study the distance from $1/y$ to $0$. For all other places $\nu$ in $S$, we fix $\alpha_{\nu}=\theta$. Hence, we have 
\begin{equation}\label{eqproduct}
|y-\infty|\prod_{p \mid g}|y-\theta|_p \prod_{q \mid c_i}|y-\theta|_q\ll |y|^{-1}\prod_{j=1}^kp_j^{-a_jm_2}=y^{-1}g^{-m_2}.    
\end{equation}
Using \eqref{eqrestr}, we obtain
\[c_2g^{m_2}+c_3g^{m_3}\leq c_2g^{m_2}+c_3g^{\frac{1}{\kappa_1}m_2} \leq \max \{|c_2|, |c_3|\}g^{\frac{1}{\kappa_1}m_2}.\]
From \eqref{eq10}
\[|y|^d\leq \Delta +cg^{\frac{1}{\kappa_1}m_2+1},\] for some constant $c$. Since $d\geq 2$ and $y\in \Z$,
\[H(y)=y\ll g^{\frac{m_2}{\kappa_1d}}\ll  g^{\frac{m_2}{2\kappa_1}}.\]
Now we will deduce that 
\begin{equation}\label{eq12}
g^{-m_2}\ll y^{-3/2}.    
\end{equation}
Since $\kappa_1>3/4$, for large values of $m_3$, we have that 
\[\kappa_1m_3>\frac{3}{2d}(m_3+1).\]
Also from \eqref{eq4} and \eqref{eqrestr}, we get
\[y^{3/2}\ll g^{\frac{3}{2d}(m_3+1)}\ll g^{m_2}.\]
Hence, we obtain \eqref{eq12}.
Combining all the above inequalities, \eqref{eqproduct} implies that 
\[|y-\infty|\prod_{p \mid g}|y-\theta|_p \prod_{q \mid c_i}|y-\theta|_q\ll H(y)^{-5/2}\]
which contradicts Theorem \ref{Roth}.

\subsection{Case-4} Our solutions run through an infinite sequence satisfying \eqref{eqrestr} and moreover such that $m_1\xrightarrow{} \infty$.

We have
\[
y^d=1+c_1g^{m_1}+c_2g^{m_2}+c_3g^{m_3}\equiv 1\bmod p^{am_1},
\] where $a=\nu_p(g)$ for some fixed
$p$ dividing $g$.
Since $(d, \varphi(p^{am_1}))=1$, we have 
\[
y\equiv 1 \bmod p^{am_1}.
\] For $m_1\geq 3$, the binomial series for $(1+c_1g^{m_1})^{1/d}$ converges in $\mathbb{Q}_p$ to an element $\zeta_p\in \mathbb{Q}_p$  and we have 
\[
\zeta_p\equiv 1\bmod p^{am_1}.
\]
Since $\zeta_p^d=1+c_1g^{m_1}\equiv y^d\bmod p^{am_2}$, we get
\begin{equation}\label{case4eq2}
    y\equiv\zeta_p\bmod p^{am_2-1}.
\end{equation}

We now choose the integers $r,s$ as the largest integers such that $r<\frac{1}{2m_1}(m_2-\frac{m_3}{d})$ and $s<\frac{1}{2m_1}(m_2+\frac{m_3}{d})$. Using \cite[ Lemma 2.5]{corvaja2013finiteness}, we may assume that either $G(-c_1g^{m_1})\neq yH(-c_1g^{m_1})$ or $G^*(-c_1g^{m_1})\neq yH^*(-c_1g^{m_1})$. Otherwise we have $G^*(-c_1g^{m_1})H(-c_1g^{m_1})-H^*(-c_1g^{m_1})G(-c_1g^{m_1})=-C(-c_1g^{m_1})$, where the left hand side equation becomes zero but the right hand side equation is non-zero. Therefore, by replacing if necessary $r,s$ respectively with $r+1,s+1$, we may assume that
\begin{equation}\label{case4eq3}
    G(-c_1g^{m_1})\neq yH(-c_1g^{m_1})
\end{equation}
and that 
\begin{equation}\label{case4eq1}
    \left| r-\frac{1}{2m_1}\left(m_2-\frac{m_3}{d}\right) \right|\leq 1, \quad \left| s-\frac{1}{2m_1}\left(m_2+\frac{m_3}{d}\right) \right|\leq 1.
\end{equation}
Let now $A$ be a positive integer as in Lemma \ref{prelem2} such that the polynomials such that $AG(z), AH(z)$ have integer coefficients bounded in absolute value by $B:=(2d)^{4t}$. Now consider the series
\[
(1-g^2z)^{1/d}=1+\frac{1}{d}(g^2z)+\frac{\left(\frac{1}{d}\right)\left(\frac{1}{d}-1\right)}{2!}(g^2z)^2+\cdots+\frac{\left(\frac{1}{d}\right)\cdots\left(\frac{1}{d}-(n-1)\right)}{n!}(g^2z)^n+\cdots .
\]
Since $(d,\phi(p^n))=1$ for every $n$, in particular $p\nmid d$. Hence $d$ is a $p$-adic unit and $|\frac{1}{d}|_p=1$. For $n\geq 1$, the coefficient of $z^n$ in the above series is
\[
\frac{1}{n!}\prod_{j=0}^{n-1}\left(\frac{1}{d}-j\right)g^{2n}.
\] Since $p\mid g$, we have
\[
\nu_p(g^{2n})=2n\nu_p(g)\geq 2n.
\]
Moreover, since $\frac{1}{d}\in\mathbb{Z}_p$, each factor $\frac{1}{d}-j$ belongs to $\mathbb{Z}_p$. Therefore,
\[
\nu_p\left(\prod_{j=0}^{n-1}\left(\frac{1}{d}-j\right)\right)\geq 0.
\]
Consequently,
\[
\nu_p\left(
\frac{1}{n!}
\prod_{j=0}^{n-1}\left(\frac{1}{d}-j\right)g^{2n}
\right)
\geq 2n-\nu_p(n!).
\]
Since
\[
\nu_p(n!)=\frac{n-S_p(n)}{p-1}
\leq \frac{n-1}{p-1}<n,
\]
where $S_p(n)\geq 1$ denotes the sum of the digits of $n$ in base $p$. Hence
\[
2n-\nu_p(n!)>n>0.
\]
Thus every coefficient of $(1-g^2z)^{1/d}$ is a $p$-adic integers and it follows from  \eqref{prelem2eq2} that $(g^2)^{t+1}AE(g^2z)$ has $p$-adic integral coefficients. 
Putting $z=-c_1g^{m_1}$ in \eqref{prelem2eq2} and recalling our notation $\zeta_p$ for the $p$-adic valuation of the binomial series for $(1+c_1g^{m_1})^{1/d}$, we obtain (using $m_1>2$)
\begin{align*}
    |AG(-c_1g^{m_1})-\zeta_pAH(-c_1g^{m_1})|_p&=|(-c_1g^{m_1})^{t+1}AE(-c_1g^{m_1})|_p\\
    &=|-c_1^{t_1+1}g^{(m_1-2)(t+1)}g^{2(t+1)}AE(-c_1g^{m_1-2}g^2)|_p\\
    &\leq p^{-a(m_1-2)(t+1)},
\end{align*} since $(c_1,g)=1$ and $|g^{2(t+1)}AE(-c_1g^{m_1-2}g^2)|_p\leq 1$, where $|\cdot |_p$ denotes the standard $p$-adic absolute value in $\mathbb{Q}_p$ and $a=\nu_p(g)$. Now, in view of \eqref{case4eq1}, we have $t+1=r+s+1\geq \frac{m_2}{m_1}-1$ and choose $k$ to be the smallest integer greater than or equal to $\frac{\ln g}{\ln p}\geq 1$ so that
\[
|AG(-c_1g^{m_1})-\zeta_pA(-c_1g^{m_1})|_p\leq p^{ak(m_1-m_2)\frac{(m_1-2)}{m_1}}.
\]
Using also \eqref{case4eq2} to eliminate $\zeta_p$ (note that $\frac{a(m_2-m_1)(m_1-2)}{m_1}\leq a(m_2-1)\leq am_2-1$) and recalling \eqref{case4eq3} we have
\begin{equation}\label{eqforp}
    0<|AG(-c_1g^{m_1})-yAH(-c_1g^{m_1})|_p\leq p^{ak(m_1-m_2)\frac{(m_1-2)}{m_1}}.
\end{equation}

Now we estimate the ordinary absolute value. By Lemma \ref{prelem2}, we have $G(z)$ and $H(z)$ are polynomials of degree $s$ and $r$ respectively.
Let 
\[
G(z)=a_0+a_1z+\cdots+a_sz^s \quad \text{and} \quad H(z)=b_0+b_1z+\cdots+b_rz^s.
\]
Then by the said bound on the coefficients of $AG(z), AH(z)$ we have
\begin{align*}
    |AG(-c_1g^{m_1})|&=|A(a_0+a_1(-c_1g^m_1)+\cdots+a_s(-c_1g_1^{m_1})^s)|\\
   & \leq B(1+c_1g^{m_1}+\cdots+c_1^sg^{m_1s})=C_1Bc_1^sg^{m_1s}.
\end{align*}
Similarly, 
\[
|AH(-c_1g^{m_1})|\leq C_2Bc_1^rg^{m_1r}.
\]
Also,
\[
y^d\leq C_3g^{m_3}\leq g^{m_3+C_4}, \quad\text{ where } C_3:=1+c_1+c_2+c_3
\]
and thus, $|y|\leq g^{\frac{m_3+C_4}{d}}$. Hence
\begin{align*}
    |AG(-c_1g^{m_1})-yAH(-c_1g^{m_1})|&\leq C_1Bc_1^sg^{m_1s}+g^{\frac{m_3+C_4}{d}}C_2Bc_1^rg^{m_1r}\\
    &=C_5Bc_1^sg^{m_1s}+C_5Bc_1^rg^{m_1r+\frac{m_3+C_4}{d}},
\end{align*}
where $C_5:=\max\{ C_1,C_2\}$. Without loss of generality assuming $s>r$, the above inequality becomes
\begin{equation*}
    |AG(-c_1g^{m_1})-yAH(-c_1g^{m_1})|\leq C_5(2d)^{4t}c_1^s\left( g^{m_1s}+g^{m_1r+\frac{m_3+C_4}{d}}\right).
\end{equation*}

On the other hand, inequalities \eqref{case4eq1} above yield
\[
t\leq \frac{m_2}{m_1}+2, \quad m_1s\leq \frac{m_2}{2}+\frac{m_3}{2d}+m_1, \quad m_1r+\frac{m_3}{d}\leq \frac{m_2}{2}+\frac{m_3}{2d}+m_1
\]
\[
s<1+\frac{1}{2m_1}\left(m_2+\frac{m_3}{d}\right).
\] Thus we get
\begin{align}\label{eqabsolute}
    |AG(-c_1g^{m_1})-yAH(-c_1g^{m_1})|&\leq  C_5(2d)^{4\left(\frac{m_2}{m_1}+2\right)}c_1^{1+\frac{1}{2m_1}\left(m_2+\frac{m_3}{d}\right)}g^{\frac{m_2}{2}+\frac{m_3}{2d}+m_1} (1+g^{C_4/2})\\
    &=C_6(2d)^{4\left(\frac{m_2}{m_1}+2\right)}c_1^{\frac{1}{2m_1}\left(m_2+\frac{m_3}{d}\right)}g^{\frac{m_2}{2}+\frac{m_3}{2d}+m_1} \notag \\
    &\leq C_6(2d)^{4\left(\frac{m_2}{m_1}+2\right)}c_1^{\frac{1}{2m_1}\left(m_2+\frac{m_3}{d}\right)}p^{ak\left(\frac{m_2}{2}+\frac{m_3}{2d}+m_1\right)}, \notag
\end{align}
where $C_6:=C_5c_1(1+g^{C_4/2})$. Since 
\[
p^{\nu_p(AG(-c_1g^{m_1})-yAH(-c_1g^{m_1}))}\leq |AG(-c_1g^{m_1})-yAH(-c_1g^{m_1})|
\] and using \eqref{eqforp} and \eqref{eqabsolute}  we obtain
\[
 C_6(2d)^{4\left(\frac{m_2}{m_1}+2\right)}c_1^{\frac{1}{2m_1}\left(m_2+\frac{m_3}{d}\right)}p^{ak\left(\frac{m_2}{2}+\frac{m_3}{2d}+m_1+(m_1-m_2)\frac{(m_1-2)}{m_1}\right)}\geq 1.
\]
That is,
\begin{equation}\label{eqgeq1}
    C_6(2d)^{4\left(\frac{m_2}{m_1}+2\right)}c_1^{\frac{1}{2m_1}\left(m_2+\frac{m_3}{d}\right)}p^{\frac{2m_2}{m_1}}p^{ak\left(\frac{-m_2}{2}+2m_1+\frac{m_3}{2d}\right)}\geq 1.
\end{equation}
However, in view of the present assumptions we have
\[
\frac{-m_2}{2}+\frac{m_3}{2d}+2m_1\leq \frac{-\kappa_1m_3}{2}+\frac{m_3}{4}+2\kappa_2m_3=-\kappa_3m_3,
\] where $\kappa_3=\frac{\kappa_1}{2}-\frac{1}{4}-2\kappa_2>0$.
Since $m_1\xrightarrow{}\infty$ along our sequence, this last inequality shows that \eqref{eqgeq1} is eventually inconsistent. This completes the proof of Theorem \ref{mainthm}.

 \section{Proof of Theorem \ref{thm2}}

Assume that
\begin{equation}\label{thm2eq1}
    y^q=1+g^{m_1}+g^{m_2}+g^{m_3}+g^{m_4},
\end{equation}
with integers $m_4>m_3>m_2>m_1>0$, $g\ge2$, $y>1$, and $q$ prime. 

{\bf Case-1} ($m_4\geq 2m_3$). 
In this case, 
\begin{align*}
    1-g^{m_4}y^{-q}&=y^{-q}(y^q-g^{m_4})=y^{-q}(g^{m_3}+g^{m_2}+g^{m_1}+1)\\
    &\leq y^{-q}(3g^{m_4/2}+1)<y^{-q}(3y^{q/2}+1).
\end{align*}We know that if $|x-1|<\frac{1}{2}, \text{ then } |\log x|<2|x-1|$. Therefore, 
\[
q\log y-m_4 \log g<y^{-q}(3y^{q/2}+1) 
\] and taking $\Lambda=|m_4\log g-q \log y|,$
\begin{equation}\label{eq3}
    \log \Lambda\ll -q\log y.
\end{equation}
Let us apply Proposition \ref{preprop1} with $b_1=q$, $b_2=m_4$, $\alpha_1=y$, and $\alpha_2=g$.
The numbers $g$ and $y$ are multiplicatively independent. Otherwise, there exists a prime $p$ dividing both $g$ and $y$ but it would lead to a contradiction in equation \eqref{thm2eq1}.

From $m_4\log g < q\log y$, we obtain
\[
b' =\frac{q}{\log g}+\frac{m_4}{\log y}=\frac{q\log y+m_4\log g}{\log g\log y}=\frac{2q}{\log g}.
\]
Applying Proposition \ref{preprop1} to inequality \eqref{eq3}, we obtain either
\[
\log\!\left(\frac{q}{\log g}\right) \ll 1,
\]
or
\[
q\log y \ll \log^{2}\!\left(\frac{q}{\log g}\right)\log y \log g.
\]
Therefore, in  both cases we conclude that
\[
q \ll \log g,
\]
where the implied constant is absolute.

{\bf Case-2} ($m_3\geq m_4/2$).

Since $m_4>m_3>m_3>m_1$, we have
\[
y^q \equiv 1 \pmod{g^{m_1}}.
\] 
First, suppose that $q$ divides $g$. From Lemma~\ref{prelem1} we obtain
\[
y^q\equiv 1\pmod{q^2}.
\] Further denoting by $q^\ell$ be the largest power of $q$ dividing $g$. This implies that 
\[
\nu_q(y^q-1)=\nu_q(g^{m_4}+g^{m_3}+g^{m_2}+g^{m_1})=\min\{m_4\nu_q(g),m_3\nu_q(g), m_2\nu_q(g), m_1\nu_q(g)\}=m_1\nu_q(g).
\]
Hence
\[
m_1\ell=\nu_q(y^q-1)\ge2,
\]
and again by Lemma~\ref{prelem1}, we have $\nu_q(y-1)=m_1\ell-1$. That is,
\[
y \equiv 1 \pmod{q^{m_1\ell-1}}.
\] 
By our assumption $q$ is coprime to $\varphi(g)$, so we have $\ell=1$. Then we derive that 
\[
y \equiv 1 \pmod{q^{m_1-1}} \text{ and } y \equiv 1 \pmod{\frac{g^{m_1}}{q}}.
\]
Also we have $$\frac{g^{m_1}}{q}\leq y-1\leq y$$ and hence
\begin{equation}\label{eq1}
    m_1\log g-\log q\leq \log y.
\end{equation}
Since $y^q \equiv 1\pmod{q^2}$, we get $2\log q\leq m_1 \log  g$. Using \eqref{eq1}, we obtain 
\begin{equation}\label{eq2}
    \frac{m_1\log g}{2}\leq \log y.
\end{equation}
Now consider $q$ does not divide $g$. Since by our assumption $q$ is coprime to $\varphi(g)$. Then we have 
\[
y \equiv 1\pmod{g^{m_1}}
\] and obtain the same inequality as in \eqref{eq2}.

Now apply Proposition \ref{preprop2} to 
\[
\Lambda= (g^{m_4}+g^{m_3}+g^{m_1}+1)-y^q
\] and $m=g^{m_1}/q$ or $m=g^{m_1}$ depending on whether $q$ divides $g$ or not. Observe that in either case, $\log m\geq \dfrac{m_1\log g}{2}$. Since $\log y \gg m_1\log g$, Proposition \ref{preprop2}  yields that
\[
\nu_m(\Lambda)
\ll
\log^2\!\left(\frac{q}{m_1\log g}\right)
\frac{\log y}{m_1^3\log^3 g}
+
\frac{\log y}{m_1\log g}.
\] On the other hand,
\[
\nu_m(\Lambda)=\nu_m(g^{m_2})=\left[ \frac{m_2}{m_1}\right]\gg \frac{\kappa q \log y}{m_1 \log g},
\] where the last inequality is obtained from the assumption that $m_2\geq \kappa m_3$ and $m_3\geq m_4/2$. Comparing the bounds of $\nu_m(\Lambda)$ yields an inequality forcing $q$
to be bounded. Consequently,
\[
q \le C(\kappa)
\]
for some effectively computable constant $C(\kappa)$ depending only on
$\kappa$. This completes the proof.
\section{Concluding remark}
In Theorem \ref{mainthm}, we proved finiteness of $(y, m_1, m_2, m_3)$ for fixed integer $d$. However, one can show that there exists a computable constant $d_1$ such that for $d>d_1$, \eqref{eq4} has no solution. Furthermore, it will be an interesting problem to remove the condition $m_2\geq \kappa m_3$ in Theorem \ref{thm2}.

\end{document}